\documentclass{amsart}
\usepackage{amsthm, amsmath}

\newtheorem{theorem}{Theorem}

\newtheorem{proposition}[theorem]{Proposition}

\numberwithin{equation}{section}
\numberwithin{figure}{section}
\numberwithin{theorem}{section}

\usepackage[utf8]{inputenc}
\newcommand{\D}{\Omega}
\newcommand{\Cn}{\mathbb{C}^n}
\newcommand{\Pro}{\mathbf{B}_{\Omega}}
\newcommand{\ProD}{\mathbf{B}_{\mathbb{D}}}
\newcommand{\dbars}{\overline{\partial}^{\ast}}
\newcommand{\dbar}{\overline{\partial}}
\allowdisplaybreaks

\title{A Survey of the $L^p$ Regularity of the Bergman Projection}
\author{Yunus E. Zeytuncu}
\address[Yunus E. Zeytuncu]{University of Michigan--Dearborn, Department of Mathematics and Statistics, 2048 Evergreen Road, Dearborn, MI 48128, USA}
\email{zeytuncu@umich.edu}

\thanks{This work is supported by NSF (DMS-1659203) and by a grant from the Simons Foundation (\#353525).}
\begin{document}

\begin{abstract}
Although the Bergman projection operator $\Pro$ is defined on $L^2(\D)$, its behavior on other $L^p(\D)$ spaces for $p\not =2$ is an active research area. We survey some of the recent results on $L^p$ estimates on the Bergman projection.    
\end{abstract}

%\dedicatory{To the memory of Nick Hanges.}

\maketitle

\section{Introduction}

Let $\D$ be a bounded domain in $\Cn$ and let $L^2(\D)$ denote the space of square-integrable functions on $\D$ with respect to the Lebesgue measure (denoted by $dV$) on $\Cn$. $L^2(\D)$ is a Hilbert space with the inner product
\begin{equation*}
\left<f,g\right>=\int_{\D}f\overline{g}dV,    
\end{equation*}
and the norm
\begin{equation*}
\left\|f\right\|^2=\int_{\D}|f|^2dV.    
\end{equation*}
The subspace of holomorphic functions in $L^2(\D)$, the Bergman space, is denoted by $A^2(\D)$. If $K$ is a compact subset of $\D$ then there exists a constant $C_K$ such that for all $f\in A^2(\D)$
\begin{equation*}
\sup_{z\in K}|f(z)|\leq C_K\int_{\D} |f|^2dV.
\end{equation*}
This so-called Bergman inequality is a consequence of the Cauchy integral formula. It follows from this inequality that $A^2(\D)$ is a closed subspace of $L^2(\D)$, and therefore an orthogonal projection operator exists from $L^2(\D)$ onto $A^2(\D)$. This projection operator is the Bergman projection $\Pro$  on $\D$. Furthermore, by the Riesz representation theorem one concludes that $\Pro$ is an integral operator, that is 
\begin{equation*}
\Pro(f)(z)=\int_{\D}B(z,w)f(w)dV,    
\end{equation*}
for $f\in L^2(\D)$ where the kernel $B(z,w)$ is called the Bergman kernel. If $\{e_n(z)\}_{n=0}^{\infty}$ is an orthonormal basis for $A^2(\D)$ then 
\begin{equation*}
B(z,w)=\sum_{n=0}^{\infty}e_n(z)\overline{e_n(w)}.    
\end{equation*}
This representation shows that $B(z,w)=\overline{B(w,z)}$ for $z,w\in \D$. $\Pro$ is a linear operator by definition, and it is self-adjoint since it is an orthogonal projection. See also \cite{KrantzBook, Jar, Krantz2} for further definitions and basic properties.

Although defined between $L^2$-spaces, the study of mapping properties of $\Pro$ on other spaces is an active research area. In particular, the behavior of $\Pro$ on Sobolev spaces $W^{k,p}(\D)$, $L^p(\D)$ spaces for $p\not =2$, and similar spaces defined by weights different than the Lebesgue measure, are studied in the literature: see \cite{ZeyTran} and the references therein. The main reason for such investigations is the connection between the mapping properties of $\Pro$ and the geometric properties of the domain $\D$. A small change in the geometry of the domain can result in significant changes in the mapping properties. For example, on the generalizations of Hartogs triangle a minor change in the definition of the domain results in a major change in the $L^p$ boundedness range of $\Pro$, see Section 3 and \cite{Edholm1, Edholm2}. Similarly, on worm domain the winding of the domain relates to the $L^p$-Sobolev boundedness range of $\Pro$, see \cite{BoasStraubeSurvey, BarSon}. Furthermore, there are many necessary and sufficient geometric conditions in the literature that relate to the estimates on $\Pro$. However, some simple cases still remain open. In particular, it is not yet known on a smooth bounded pseudoconvex domain whether a finite type assumption is sufficient for $L^p$ regularity, see the list of similar open problems in the last section.

Before we focus on the $L^p$ regularity of $\Pro$, we highlight an analogous (and more complete) study on Sobolev spaces. Indeed, the behavior of $\Pro$ on Sobolev spaces (especially $L^2$-Sobolev spaces $W^k(\D)$) is an interesting section of the regularity studies. One reason for this is that the Sobolev regularity of $\Pro$ is closely related to the Sobolev regularity of the $\dbar$-Neumann operator $N$. Indeed, the Kohn formula (see \cite{Kohn}) gives the following representation 
\begin{equation*}
    \Pro (f)=f-\overline{\partial}^{\ast}N\dbar f.
\end{equation*}
for any $f\in L^2(\D)$. We refer to \cite{BoasStraubeSurvey} for the definitions of the operators above and a comprehensive account of Sobolev regularity. We only highlight one important aspect here. If $\Pro$ maps the Sobolev space $W^s(\D)$ onto itself for all $s>0$ we say that $\Pro$ is exactly regular. If $\Pro$ maps $C^{\infty}(\overline{D})$ to itself we say that $\Pro$ is globally regular (or satisfies Condition-R). On a smooth bounded domain $\D$, exact regularity implies global regularity by the Sobolev embedding theorem. Furthermore, exact regularity turns out to be a sufficient condition for the smooth extension of biholomorphisms. Indeed, it was shown in \cite{Bell81} (see also \cite{BellLig} for a precursor) that if $F$ is a biholomorphism between two bounded domains with smooth boundary $$F:\D_1\to \D_2$$ and if the Bergman projection of one of the domains is globally regular then $F$ extends to a diffeomorphism between the closures of these two domains. An earlier version of this result on strongly pseudoconvex domains was obtained by Fefferman in \cite{Fefferman}. This extension is particularly important in the investigation of the class of biholomorphically equivalent domains in $\Cn$. We again refer the reader to \cite{BoasStraubeSurvey} for the details. 

Now we shift our attention back to the $L^p$ regularity of $\Pro$. The main question is to determine the values of $p$ such that $\Pro$ extends to a bounded operator  on $L^p(\D)$. By definition $\Pro$ maps $L^2(\D)$ into itself. If $\Pro$ maps $L^p(\D)$ into itself for some $p>2$ then it is also bounded on $L^q(\D)$ for the conjugate exponent $q$: $\frac{1}{p}+\frac{1}{q}=1$ by the duality of $L^p(\D)$-spaces and self-adjointness of $\Pro$. Furthermore, $\Pro$ is linear and by interpolation one can also conclude that if $\Pro$ is bounded on $L^p(\D)$ onto itself for some $p>2$ then it is also bounded on $L^r(\D)$ for all $q\leq r \leq p$. In other words, for any given domain $\D$, we can associate an interval $I_{\D}$ that at least contains $\{2\}$ and contained in $(1,\infty)$ such that $\Pro$ is bounded on $L^p(\D)$ if and only if $p\in I_{\D}$. Calculating $I_{\D}$ for a specific $\D$ is non-trivial. It can be as small as the singleton $\{2\}$ (see Theorem \ref{Zey1}), as big as the full interval $(1,\infty)$ (see Section 2), and it can be any conjugate symmetric sub-interval that contains $2$ (see Theorem \eqref{p_0}).

We note that the boundedness of $\Pro$ from $L^p(\D)$ to $L^p(\D)$ is the best estimate one can hope in this context, since $\Pro$ equals the identity on the space of holomorphic functions in $L^p(\D)$. We also note that on domains (with smooth enough boundary) where exact regularity holds, by using the Sobolev embedding theorem, one can obtain a crude estimate of the form: for any $p>2$, $\Pro$ maps $L^p(\D)$ to $L^{p-\delta}(\D)$ for some $\delta>0$ that depends on $p$ and the dimension of the ambient space. We refer to \cite{HarZey} for further discussion on how to calculate a small $\delta$ when $p$ is sufficiently close to $2$.

In the rest of this note, we go over some recent results in the literature on this problem. We also list some open problems (see the last section). In our discussion here we only consider the Lebesgue measure on $\D$ without any extra weight factor. The study of weighted (or unweighted) Bergman projection on weighted $L^p$ or Sobolev spaces is also an active direction. We leave a survey of those results to another work.

\section{Boundedness}
The earliest $L^p$ estimates on the Bergman projection were obtained on the unit disc $\mathbb{D}$ in $\mathbb{C}^1$, see \cite{Zah} for the original proof, and see \cite{KeheBook} for a detailed account on the unit disc. In this case the Bergman kernel is explicitly known. Indeed, the set of  monomials $\left\{\sqrt{\frac{n+1}{\pi}}z^n\right\}_{n=0}^{\infty}$ is an orthonormal basis for $A^2(\mathbb{D})$ and 
\begin{equation*}
B_{\mathbb{D}}(z,w)=\sum_{n=0}^{\infty}\frac{n+1}{\pi}z^n\overline{w}^n=\frac{1}{\pi}\frac{1}{(1-z\overline{w})^2}.    
\end{equation*}
Using this representation and Schur's lemma, one can show that $(1,\infty)\subset I_{\mathbb{D}}$. Furthermore, one can also show explicitly that $\ProD$ is unbounded on $L^1(\mathbb{D})$ and $L^{\infty}(\mathbb{D})$. In other words, one concludes that $I_{\mathbb{D}}=(1,\infty)$. See also \cite{Zhao} for a detailed account with proofs--not only for $n=1$, but also for the unit ball in $\Cn$ and $n>1$.

There are many other planar domains $\D\subset \mathbb{C}$ on which the Bergman projection is bounded on $L^p(\D)$ for all $p\in(1,\infty)$. In many cases of simply connected domains one can prove such a result by using the Riemann mapping theorem. Indeed, if two domains are conformally equivalent then the Bergman kernels can be related by the conformal map.
If $F:D_{1}\to D_{2}$ is a conformal map then the Bergman kernels are related by 
\begin{equation}\label{kernel}
    B_{D_1}(z,w)=F'(z)B_{D_2}(F(z),F(w))\overline{F'(w)}.
\end{equation}
In particular, if the second domain $D_2$ is the unit disc $\mathbb{D}$ then the Bergman kernel of $D_1$ is given by
\begin{equation}
    B_{D_1}(z,w)=\frac{F'(z)\overline{F'(w)}}{\pi (1-F(z)\overline{F(w)})^2}.
\end{equation}
Furthermore, \eqref{kernel} implies the following relation between the operators (we denote the inverse of the conformal map $F$ by $G$), for $\phi \in L^p(D_1)$
\begin{equation*}
G'(w)\mathbf{B}_{D_1}(\phi)(G(w))=\mathbf{B}_{\mathbb{D}}(\phi(G)G')(w).
\end{equation*}
Using this representation, one can conclude that the Bergman projection $\mathbf{B}_{D_1}$ is bounded on $L^p(D_1)$ for some $p\geq 1$ if and only if the Bergman projection $\mathbf{B}_{\mathbb{D}}$ on the unit disc is bounded on the weighted space $L^p(\mathbb{D},|G'|^{2-p})$. 
Indeed, for $\phi \in L^p(D_1)$, we have $\phi(G)G'\in L^p(\mathbb{D},|G'|^{2-p})$ with the equality
\[
||\phi||_{L^p(D_1)}^p=\|\phi(G)G' \|_{L^p(\mathbb{D},|G'|^{2-p})}^p
\]
and
\begin{align*}
\| \mathbf{B}_{D_1}(\phi) \|^p&=\int_{D_1}|\mathbf{B}_{D_1}\phi (z)|^pdA(z)\\
&=\int_{\mathbb{D}}|\mathbf{B}_{D_1}\phi (G(w))|^p|G'|^2dA(w)\\
&=\int_{\mathbb{D}}|G'(w)\mathbf{B}_{D_1}\phi (G(w))|^p|G'|^{2-p}dA(w)\\
&=\| \ProD \left(\phi(G)G'\right)\|^p_{|G'|^{2-p}}
\end{align*}
On the other hand, the following was observed in \cite{Bekolle, LanSte04}. 
\begin{theorem}\label{Ap}
The Bergman projection $\mathbf{B}_{\mathbb{D}}$ on the unit disc is bounded on a weighted space $L^p(\mathbb{D},\omega)$ for some weight $\omega$ if and only if $\omega\in \mathcal{A}_p^+$. 
\end{theorem}
The class $\mathcal{A}_p^+$ is a version of the Muckenhoupt class \cite{Muck} for the Riesz transform on the real line. A weight $\sigma$ on $\mathbb{D}$ is a non-negative locally continuous function. We say $\sigma$ is in $\mathcal{A}_p^+$, if there exists $C>0$ such that for any Carleson tent $T_z$ for $z\in \mathbb{D}$
\begin{equation}\label{BB}
\frac{1}{|\mathbb{D}\cap T_z|^p}\left(\int_{\mathbb{D}\cap T_z}\sigma dA\right)\left(\int_{\mathbb{D}\cap T_z}\sigma^{\frac{1}{1-p}} dA\right)^{p-1}\leq C     
\end{equation}
where $|\mathbb{D}\cap T_z|$ denotes the Lebesgue measure of the set. A Carleson tent $T_z$ for $z\not=0$ is defined by
\begin{equation*}
T_z=\left\{ w\in \mathbb{D}~:~ \left|1-\overline{w}\frac{z}{|z|}\right|<1-|z|\right\}    
\end{equation*}
and when $z=0$ the tent is the whole $\mathbb{D}$. Carleson tents can be also defined by using the Bergman metric. The Bergman kernel satisfies nice off-diagonal estimates on these tents. Therefore, the mapping properties can be controlled by controlling the weights on these sets. The estimates \eqref{BB} are sometimes referred as the Bekoll\'e-Bonami estimates. We also refer to \cite{HuoW1, HuoW2} for a generalization on pseudoconvex domains.

This chain of equivalences implies that smoothness of the conformal map $F$ on the closure of the domains determines the regularity of the Bergman kernel and therefore also the estimates on the projection operator. In particular, if $F: D_1\to \mathbb{D}$ is the Riemann mapping function with the inverse map denoted by $G$, then the Bergman projection $\mathbf{B}_{D_1}$ is bounded on $L^p(D_1)$ for some $p>2$ if and only if $|G'|^{2-p}\in \mathcal{A}_p^+$.
In this sense, $L^p$ estimates on planar domains overlap with geometric function theory and the regularity of the conformal maps. These connections lead to the following general non-degeneracy statement in \cite{Hed} on planar domains. 

\begin{theorem}\cite[Theorem 3.1]{Hed}
 If $\D$ is a proper simply connected planar domain, then there exists $p_0>2$ such that $(q_0,p_0)\subset I_{\D}$. In other words, the Bergman projection is bounded on some $L^p$ spaces for $p\not =2$.
\end{theorem}
The proof follows from the observation that for any such $\D$ the Riemann mapping function to a small power belongs to some $\mathcal{A}_p^+$ class for $p$ close enough to $2$. We will see in the next section that such low-level $L^p$ regularity is a phenomenon on planar domains, and it fails in several variables. 
 
In higher dimensions, if the domain has nice geometric structure then the $L^p$ mapping interval again will be as big as $(1,\infty)$. Indeed, Forelli and Rudin \cite{ForelliRudin} extended the planar results to the ball and polydisc in higher dimensions. In the direction of more general domains with less symmetry, Phong and Stein \cite{PhongStein}, showed that if $\D$ is a bounded, smooth, strongly pseudoconvex domain then $I_{\D}=(1,\infty)$. Their approach relies on kernel estimates and the general theory of singular integral operators. This approach later resulted in similar conclusions on more general domains. More specifically, if $\D$ satisfies one of the geometric conditions:
\begin{itemize}
\item smooth, finite type in $\mathbb{C}^2$ \cite{McNealC2},
\item smooth, convex, finite type in $\mathbb{C}^n$ \cite{McnSteConvexFinite},
\item smooth, decoupled in $\Cn$ \cite{McNealDecoupled},
\item smooth, finite type, and diagonalizable Levi form in $\Cn$ \cite{Charpentier}
\item $C^2$-smooth, strictly pseudoconvex \cite{LanSte12},
\end{itemize}{}
then $I_{\D}=(1,\infty)$. In these results, one of the key ingredients is the off-diagonal pointwise estimates on the Bergman kernel in terms of the distance to the boundary. Once these technical estimates are obtained then other tools of harmonic analysis are employed. We also refer to \cite{McNealSingular} for an account of these results from the theory of singular integral operators. We want to point out that although the results above cover a large family of domains, it is still not known yet on a smooth bounded pseudoconvex domain if finite type assumption is sufficient for $L^p$ regularity, see the open problem $\#3$ in the last section.

\section{Unboundedness}

In most of the cases in the previous section the $L^p$ mapping interval was $(1,\infty)$. There are examples of domains on which this interval is strictly smaller. Such examples can be created in the plane by taking advantage of the Riemann mapping theorem. As we mentioned above the $L^p$ regularity interval depends on the smoothness of the Riemann mapping function on the closure. This approach requires manipulating the Riemann mapping function by making the domain more singular. As indicated above by Hedenmalm's theorem, there is a limit how much one can shrink $I_{\D}$ on planar domains.

On the other hand, one can find domains in higher dimensions on which the $L^p$ boundedness interval is just $\{2\}$, the smallest possible. First such example appears in \cite{Barrett} where the domain is smooth but not pseudoconvex. In order to define this domain, one needs to start with the following smooth functions $r_1, r_2$, and $c$ on $[1,6]$ with the properties
\begin{align*}
r_1(x)&=
\begin{cases}
3-\sqrt{x-1}& \text{ near }x=1, \text{ and decreasing on } [1,2]\\
1& \text{ on }[2,5] \text{ and increasing on }[5,6]\\
3-\sqrt{6-x}& \text{ near }x=6;
\end{cases}\\
r_2(x)&=
\begin{cases}
3+\sqrt{x-1}& \text{ near }x=1, \text{ and increasing on } [1,2]\\
4& \text{ on }[2,5] \text{ and decreasing on }[5,6]\\
3+\sqrt{6-x}& \text{ near }x=6;
\end{cases}\\
c(x)&=
\begin{cases}
0& \text{ on }[1,2], \text{ and decreasing on } [2,3]\\
\exp\left({-|x-3|^{-1/2}}\right)-1& \text{ near }x=3 \text{ and increasing on }[3,4]\\
-\exp\left({-|x-4|^{-1/2}}\right)+1& \text{ near }x=4 \text{ and decreasing on }[4,5]\\
0& \text{ on }[5,6].
\end{cases}
\end{align*}{}
Now, let
\begin{equation*}
    \D=\left\{(z_1,z_2)\in \mathbb{C}^2:~1<|z_2|<6, |z_1|<r_2(|z_2|), |z_1-c(|z_2|)|>r_1(|z_2|)\right\}.
\end{equation*}
It is a smooth domain where the fibers are annuli with moving centers. On $\D$, one can show that for any $p>2$ there exists a compactly supported smooth function $\phi_p$ on $\D$ such that $\Pro(\phi_p)\not \in L^p(\D)$.

Later, in \cite{ZeyTran} another example is presented, that is pseudoconvex but not smooth. Indeed, for any $A\geq 0$, $B>0$, and $\alpha>0$, let
\begin{equation*}
\D_{A,B,\alpha}=\left\{(z_1,z_2)\in \mathbb{C}^2:~|z_1|<1, |z_2|<(1-|z_1|^2)^A\exp\left(\frac{-B}{(1-|z_1|^2)^{\alpha}}\right) \right\}.    
\end{equation*}
\begin{theorem}\label{Zey1}
The Bergman projection $\mathbf{B}_{\D_{A,B,\alpha}}$ is bounded on $L^p(\D_{A,B,\alpha})$ if and only if $p=2$. 
\end{theorem}
The proof uses the orthogonal decomposition of $A^2(\D_{A,B,\alpha})$ into weighted Bergman spaces on the unit disc. Then one can show that the corresponding weighted Bergman projections are bounded on the corresponding weighted $L^p$ spaces if and only if $p=2$. Similar weighted results appear in \cite{Dostanic, ZeyTran}. We note that in both examples the exponential decay plays an essential role in the proofs.

Similar non-smooth but pseudoconvex examples later appeared in \cite{Edholm2, Edholm1}. These domains are quite similar to the Hartogs triangle in $\mathbb{C}^2$. However, this small change in the definition has resulted in drastic changes in the $L^p$ mapping properties. Indeed, for an irrational $\gamma>0$, let
\begin{equation*}
\D_{\gamma}=  \left\{(z_1,z_2)\in \mathbb{C}^2:~0<|z_1|<1, |z_2|<|z_1|^{\gamma} \right\}.  
\end{equation*}
It turns out that on these \textit{irrational} Hartogs triangles $I_{\D_{\gamma}}=\{2\}$. We note that on the standard Hartogs triangle (that is $\gamma=1$ in the definition above), the boundedness interval is exactly $(\frac{4}{3},4)$, see \cite{ChakZey}.

Both $\D_{A,B,\alpha}$ and $\D_{\gamma}$ are Reinhardt domains and subsets of the monomials $\{z^{\alpha}\}$ for multi-indices $\alpha\in\mathbb{Z}^2$ form an orthogonal bases for the corresponding Bergman spaces. Both results above use these explicit basis elements in the proofs. For a continuous function $h:[0,1]\to [0,\infty)$ such that $h(0)=0$, the Reinhardt domain $\Omega_h$ defined by
\begin{equation}\label{Omegah}
\D_h=\left\{(z_1,z_2)\in \mathbb{C}^2:~0<|z_1|<1, |z_2|<h(|z_1|) \right\}
\end{equation}
is a natural generalization of the domains above. One can investigate the $L^p$ mapping properties of the Bergman projection on these domains and try to relate $I_{\D_h}$ to the properties of $h$. Although one may conjecture that the order of vanishing of $h(t)$ as $t\to 0^+$ is the key property, the \textit{irrational} Hartogs triangles above indicate that not only the order of vanishing but the algebraic properties of $h(t)$ at $t=0$ might be relevant in such an investigation.

It is not yet clear if such a degeneracy (that is $I_{\D}=\{2\}$) can hold on smooth pseudoconvex domains. In other words, is it true that for any bounded smooth pseudoconvex domain $\D$, there exists $p_0>2$ such that $\Pro$ is bounded on $L^p(\D)$ for any $q_0\leq p\leq p_0$?
We note that for $L^2$-Sobolev estimates such a low-level regularity holds. Indeed, for any given smooth bounded pseudoconvex domain $\D$ there exists $\varepsilon>0$ (depends on $\D$) such that $\Pro$ is bounded on $W^{\delta}(\D)$ for all $0\leq \delta<\varepsilon$, see \cite{Kohn}. In fact, this low-level Sobolev regularity is not only true for the Bergman projection, but also for the $\dbar$-Neumann operator $N$ and the canonical solution operators $\dbar N$ and $\dbars N$.

In this direction, in \cite{HarZey} the authors presented some low-level $L^p$ regularity on a class of domains for the canonical operators (but not the Bergman projection). 
\begin{theorem}\label{HarZey}
 Let $\Omega\subset\mathbb{C}^n$ be a bounded domain that is locally a transverse intersection of smooth finite type domains.  Then there exists some $p_0>2$ such that the $\dbar$-Neumann operators
  \begin{align*}
    N_q:&L^{p/(p-1)}_{(0,q)}(\Omega)\rightarrow L^p_{(0,q)}(\Omega)\\
    \dbar^* N_{q}:&L^2_{(0,q)}(\Omega)\rightarrow L^p_{(0,q-1)}(\Omega)\\
    \dbar N_{q-1}:&L^2_{(0,q-1)}(\Omega)\rightarrow L^p_{(0,q)}(\Omega)
  \end{align*}
  are continuous for all $1\leq q\leq n$ and $2<p<p_0$.
\end{theorem}
The proof is a combination of estimates with good weights, in particular gain in terms of weights due to subellipticity of $N$ on such domains, and manipulating between weighted Sobolev spaces and unweighted $L^p$ spaces.

We also note more examples of domains where the Bergman projection or the canonical operators exhibit unboundedness were studied in \cite{BarSon, PelKra1, PelKra2}.
In a sense, the degenerate examples above are due to unboundedness. Since the operator is bounded on $L^2(\D)$ by definition, in order to prove something like above one needs to find appropriate test functions and explicitly compute the projections of these functions. On the other hand, when the $L^p$ boundedness range is strictly between $\{2\}$ and $(1, \infty)$ then one needs argue in two steps: first one needs to prove unboundedness of the operator for certain values of $p$, and then prove boundedness of it for certain other values of $p$. The latter usually requires off-diagonal estimates on the Bergman kernel. In this direction, the following was obtained in \cite{ZeyTran}. 
\begin{theorem}\label{p_0}
For any given $p_0>2$, let
\begin{equation*}
\D_{p_0}=    \left\{(z_1,z_2)\in \mathbb{C}^2:|z_1|<1, |z_2|<|z_1-1|^{\frac{2}{p_0-2}} \right\}. 
\end{equation*}
Then the $L^p$ regularity interval $I_{\D_{p_0}}$ is exactly this pre-determined interval $(q_0,p_0)$, where $\frac{1}{q_0}+\frac{1}{p_0}=1$.
\end{theorem}
These domains are Hartogs domains, and one can express the Bergman space $A^2(\D_{p_0})$ as an orthogonal sum of weighted Bergman spaces on the unit disc. Indeed, it turns out that $A^2(\D_{p_0})$ is isometrically isomorphic to
\begin{equation*}
\bigoplus_{m=0}^{\infty} A^2(\mathbb{D}, |z-1|^{\frac{4(m+1)}{p_0-2}}).
\end{equation*}
One can compute the weighted Bergman kernels of $A^2(\mathbb{D}, |z-1|^{\frac{4(m+1)}{p_0-2}})$ since the weight is the norm square of a holomorphic function on $\mathbb{D}$. 
\begin{proposition}
Let $g(z)$ be a non-vanishing holomorphic function on $\mathbb{D}$. The weighted Bergman kernel $B_g(z,w)$ of the space $A^2(\mathbb{D},|g|^2)$ is given by
\begin{equation*}
B_g(z,w)=\frac{1}{g(z)}B_{\mathbb{D}}(z,w)\frac{1}{\overline{g(w)}},
\end{equation*}
where $B_{\mathbb{D}}(z,w)$ is the ordinary (unweighted) Bergman kernel on the unit disc.
\end{proposition}

Indeed, since $\frac{1}{g(z)}$ is holomorphic on the disc the set $\{\frac{z^n}{g(z)}\}_{n=0}^{\infty}$ forms an orthonormal basis for the weighted Bergman space $A^2(\mathbb{D},|g|^2)$. Using this representation and by keeping track of the $\mathcal{A}_p^+$-class membership of the weight (see Theorem \ref{Ap}) one can calculate the $L^p$ boundedness interval of the weighted Bergman projections of $A^2(\mathbb{D}, |z-1|^{\frac{4(m+1)}{p_0-2}})$. It turns out that if one of these weighted Bergman projections is unbounded on the unit disc then the operator on the Hartogs domains fails to be bounded. On the other hand, one can prove that if all the weighted projections on the disc are bounded on the respective weighted spaces $L^p(\mathbb{D}, |z-1|^{\frac{4(m+1)}{p_0-2}})$ for some $p$, then the operator $\mathbf{B}_{\D_{p_0}}$ is bounded on $L^p(\D_{p_0})$ for the same value of $p$. These two arguments lead to a proof of Theorem \ref{Ap}. In other words, any conjugate symmetric sub-interval of $(1,\infty)$ can be the $L^p$ boundedness range of a pseudoconvex domain. Later similar examples appeared in \cite{ChakZey, Edholm2, Edholm1, Chen, Huo, Beb} where the $L^p$ regularity of the domain is exactly calculated. 

We note another common theme among these results. Except the degenerate case where $I_{\D}=\{2\}$, $I_{\D}$ is always an open interval in the examples mentioned above. This observation leads to the question of whether or not there is an openness statement in this context. In other words, is it true that $I_{\D}$ is always open except for the degenerate case $I_{\D}=\{2\}$? One way of approaching such an investigation is actually proving that if $\Pro$ is bounded on some $L^p(\D)$ for some $p>2$ then it is also bounded on $L^{r}$ for some $r>p>2$. Such an openness statement holds in the context of $A^p$-weights as studied in \cite{Muck}. Indeed, if a weight $\omega$ on the unit disc is in $A^p$-class for some $p>2$ then it is also in some $A^r$-class for some $r>p$.

\section{Related Operators}
In addition to the Bergman projection, one can study the similar $L^p$ mapping properties for the related operators such as various solution operators for $\dbar$, Szeg\"o projection and Toeplitz operators. In this section, we list a few key references for the reader.

\begin{itemize}
\item For the Szeg\"o projection, statements similar to Theorems \ref{Zey1} and \ref{p_0} appeared in \cite{MunZey1, MunZey2}. Furthermore, a study of $L^p$ mapping properties on worm domains appeared in \cite{Mon1}.

\item The literature on $L^p$ estimates for solution operators is quite rich because many different solution operators exist. \cite{Sibony}, \cite{Lieb}, and \cite{Kerzman} serve as starting points. More recent results can be seen in \cite{Lee, ChenMcN}.

\item For a bounded function $\varphi$ on a domain $\D\subset \Cn$, one can define the Bergman Toeplitz operator $T_{\varphi}$ as
\begin{align*}
T_{\varphi}&:A^2(\D)\to A^2(\D),\\    
T_{\varphi}(g)&=\Pro(\varphi g).
\end{align*}
The mapping properties of $T_{\varphi}$ depend on the symbol $\varphi$ and the operator $\Pro$. \cite{CucMcN}, \cite{Abate}, and \cite{Khanh} contain the most recent $L^p$ estimates on Toeplitz operators.

\end{itemize}

\section{Open Problems}
We mentioned a few open problems and possible directions of research in the text above. We list those and a few other specific problems in this section.
\begin{enumerate}

\label{OP1}\item The classification of the weights that satisfy the Bekoll\'e-Bonami estimates \eqref{BB} on general domains is an open problem: classify weights $\sigma$ on general domains $\D$ such that the ordinary Bergman projection $\Pro$ extends to be a bounded operator on the weighted space $L^p(\D,\sigma)$ for $p>1$. Furthermore, relate the operator norm of $\Pro$ to the weight $\sigma$. Recently, \cite{HuoW1} and \cite{HuoW2} answered this question on some pseudoconvex domains on which sharp off-diagonal estimates on the Bergman kernel are known. They use a careful construction of dyadic decomposition on these domains by generalizing Carleson tents on the unit disc.

\label{OP2}\item Is it true that on any smooth pseudoconvex bounded domain $\D\subset \Cn$, there exists $\varepsilon>0$ (depends on $\D$) such that $(\frac{2+\varepsilon}{1+\varepsilon},2+\varepsilon)\subset I_{\D}$? In other words, is there a smooth pseudoconvex domain such that the $L^p$ mapping interval $I_{\D}$ degenerates to $\{2\}$?

\label{OP3}\item Characterize the smooth pseudoconvex domains $\Omega$ for which $I_{\D}=(1,\infty)$. In particular, is it true for all bounded pseudoconvex finite type domains in $\Cn$ ($n>2$) the $L^p$ boundedness range is $(1,\infty)$?

\label{OP4}\item For a given domain $\D$, except the degenerate case of $\{2\}$, can the $L^p$ boundedness interval $I_{\D}$ be a closed interval?

\label{OP5}\item Let $\Omega$ be a complete Reinhardt domain. Is the $L^p$ mapping range $I_{\Omega}$ the full interval $(1,\infty)$? Using symmetry, Boas-Straube proved that  $\Pro$ is exactly regular \cite{BoasStraubeSurvey}. See also \cite{Huo}.

\label{OP6}\item Most of the boundedness $L^p$ estimates on $\Pro$ stem from precise kernel estimates on $B_{\D}(z,w)$. Is it possible to deduce the $L^p$ boundedness estimates on $\Pro$ by using similar estimates on the canonical operators as in Theorem \ref{HarZey}? In particular, can one extract any $L^p$ estimates on smooth pseudoconvex domains on which the $\dbar$-Neumann operator is compact? 

\label{OP7}\item For the domains $\Omega_h$ (defined as in \eqref{Omegah}), investigate the relation between the $L^p$ mapping properties of $\mathbf{B}_{\D_h}$ and the properties of $h(t)$. Is it true that $I_{\D_h}$ is controlled by the behavior of $h(t)$ around $t=0$?

\end{enumerate}

\section*{Acknowledgements}
I thank Mehmet \c{C}elik, John D'Angelo, Nordine Mir, and S\"onmez \c{S}ahuto\u{g}lu %, and the anonymous referee 
for many useful comments on preliminary versions of this paper.

\bibliographystyle{alpha}
\bibliography{Survey}

\end{document}